\documentclass[11pt,twoside,a4paper]{article}

\usepackage[utf8]{inputenc}
\usepackage{amsmath,amsfonts,amssymb,amsthm}
\usepackage{enumerate,lineno,setspace,float}
\usepackage[numbers,sort&compress]{natbib}
\usepackage[small,compact]{titlesec}
\usepackage[capitalize]{cleveref} 
\usepackage{graphicx,float,url}
\usepackage{geometry}
\usepackage{fancyhdr}
\usepackage{latexsym}
\usepackage{setspace}
\usepackage{tikz}
\usepackage{caption}
\usepackage{xcolor}

\usepackage[capitalize]{cleveref}

\newtheorem{theorem}{Theorem}[section]
\newtheorem{lemma}[theorem]{Lemma}

\newtheorem{example}[theorem]{Example}

\newtheorem{observation}[theorem]{Observation}

\input amssym.def
\input amssym.tex

\usepackage{xcolor}
\usepackage[normalem]{ulem}

\newcommand{\be}{\begin{equation}}
\newcommand{\ee}{\end{equation}}
\newcommand{\ben}{\begin{equation*}}
\newcommand{\een}{\end{equation*}}
\newcommand{\bea}{\begin{eqnarray}}
\newcommand{\eea}{\end{eqnarray}}
\newcommand{\bean}{\begin{eqnarray*}}
\newcommand{\eean}{\end{eqnarray*}}
\def\non{\nonumber}

\def\qed{\hfill$\Box$\vspace{11pt}}

\def\diam{{\rm diam}}

\def\span{{\rm span}}
\def\rn{{\rm rn}}

\def\({\left(}
\def\){\right)}
\def\[{\left[}
\def\]{\right]}

\def\an{{\rm an}}

\def\lcm{{\rm lcm}}

\begin{document}

\title{Radio antipodal number for the Cartesian product of a complete graph and a path}

\author{\textbf{Payal Vasoya} \\ Gujarat Technological University, Ahmedabad - 382 424, Gujarat (India) \\ E-mail : prvasoya92@gmail.com \\ \\
\textbf{Devsi Bantva}\footnote{Corresponding author.} \\ Lukhdhirji Engineering College, Morvi - 363 642, Gujarat (India) \\ E-mail : devsi.bantva@gmail.com}

\pagestyle{myheadings}
\markboth{\centerline{Payal Vasoya and Devsi Bantva}}{\centerline{Radio antipodal number for the Cartesian product of a complete graph and a path}}

\date{}
\openup 0.8\jot
\maketitle

\begin{abstract}
A radio antipodal labeling of a graph $G$ is a function $f$ that assigns a non negative integer to each vertex in such a way that for any pair of distinct vertices $x$ and $y$, the condition $|f(x)-f(y)| \geq d-d(x,y)$ is satisfied, where $d$ and $d(x,y)$ represent the diameter of graph $G$, and distance between $x$ and $y$ in graph $G$, respectively. The span of radio antipodal labeling $f$ is defined as $\span(f) = \max\{|f(x)-f(y)| : x,y \in V(G)\}$. The radio antipodal number of graph $G$, denoted by $\an(G)$, is defined as $\an(G) := \min_{f}\{\span(f)\}$ with minimum taken over all radio antipodal labelings $f$ of $G$. A radio antipodal labeling $f$ of $G$ is optimal if $\span(f) = \an(G)$. In this paper, we give the exact radio antipodal number of the Cartesian product of a complete graph and a path. 

\end{abstract}

\section{Introduction}\label{intro}
In wireless communication networks, assigning frequencies with a minimum span to transmitters in a way that avoids interference is an important problem, known as the \textit{channel assignment problem}. This problem can be modeled using a graph model, where each transmitter is represented by a vertex, and interference between two transmitters is represented by an edge. Different frequencies (or channels) must be assigned to nearby transmitters based on their proximity to prevent interference. In 1980, Hale proposed various graph labelings in \cite{Hale} to deal the channel assignment problem. In 1992, Griggs and Yeh\cite{Griggs} introduced \textit{$L(2,1)$-labeling}, also known as \textit{distance-two labeling}. In $L(2,1)$-labeling, if two vertices are adjacent, then the labels assigned to them must differ by at least $2$; and if they are at distance $2$, then their labels must differ by at least $1$. The difference between the maximum and minimum labels assigned by an $L(2,1)$-labeling $f$ is called the span of $f$, denoted by $\span(f)$. The minimum span taken over all $L(2,1)$-labelings $f$ of $G$ is called the $\lambda$-number of graph $G$, denoted by $\lambda(G)$. That is, $\lambda(G) = \min\{\span(f) : f \text{ is an $L(2,1)$-labeling of $G$}\}$. This labeling ensures sufficient separation between labels to reduce interference among nearby transmitters.

Later, in 2001, Chartrand \emph{et al.} \cite{Chartrand} generalized this concept by introducing \textit{radio $k$-labeling}, inspired by more realistic communication constraints, in which interference can occur between any pair of vertices, not just adjacent ones. A radio $k$-labeling is a function $f : V(G) \rightarrow \{0,1,2,\ldots\}$ such that the following holds for every pair of distinct vertices $u$ and $v$ in the graph $G$:
$$
|f(u) - f(v)| \geq k + 1 - d(u,v),
$$
where $d(u,v)$ is the distance between $u$ and $v$, and $k$ is a positive integer. The span of $f$ is defined as $\span(f) = \max\{|f(u)-f(v)| : u,v \in V(G)\}$. The radio $k$-number of $G$, denoted by $\rn_k(G)$, is defined as
$$\rn_k(G) := \min_{f}\{\span(f)\},$$
where minimum taken over all radio $k$-labelings $f$ of $G$. A radio $k$-labeling $f$ of $G$ is optimal if $\span(f) = \rn_k(G)$.

Interestingly, the radio $k$-labeling problem generalizes many graph labeling problems. For $k=1$, radio $1$-labeling is a famous proper vertex coloring and $\rn_1(G) = \chi(G)$, the chromatic number of $G$. For $k=2$, radio $2$-labeling is equivalent to $L(2,1)$-labeling (or distance two labeling) and $\rn_2(G) = \lambda(G)$. For $k=d-1$ (where $d \geq 2$ is the diameter of graph $G$), radio $(d-1)$-labeling is known as radio antipodal labeling and $\rn_{(d-1)}(G) = \an(G)$, the radio antipodal number of $G$. Finally, for $k = d$, radio $d$-labeling is called radio labeling, and $\rn_d(G) = \rn(G)$, the radio number of $G$. 

A radio antipodal labeling assigns distinct labels to each pair of distinct vertices, except for pairs of antipodal vertices, which may have the same label. Therefore, a radio antipodal labeling is slightly more relaxed than a radio labeling of graphs, yet more constrained than radio $k$-labelings for $k < d$. Without loss of generality, we may assume that a radio antipodal labeling $f$ always assigns the label $0$ to some vertex, and in this case, the span of $f$ is the maximum label used by $f$. A radio antipodal labeling $f$ induced an ordering $z_1,z_2,\ldots,z_n$ of the vertices of $V(G)$ (where $n=|V(G)|$) such that 
$$0 = f(z_1) \leq f(z_2) \leq \ldots \leq f(z_n) = \span(f).$$

The study of radio antipodal labeling has evolved significantly over the years, with key contributions from various authors. Determining the radio antipodal number of paths, however, was interesting and challenging. In \cite{Zhang2}, Chartrand \emph{et al.} proved the following result about $\an(P_n)$: 
$$\an(P_n) \leq \binom{n-1}{2}+1.$$
Further, they conjectured that the value on the right-hand side of above inequality is the exact radio antipodal number of graphs. In \cite{Zhang3}, Chartrand \emph{et al.} showed that the conjecture is false and they gave another sharp upper bound for the antipodal number of $P_n$ of odd order $n \geq 7$. 
\begin{theorem}\cite{Zhang3}\label{Pn:upper} For the path $P_n$ of odd order $n \geq 7$,
$$\an(P_n) \leq \binom{n-1}{2}-\frac{n-1}{2}+4.$$
\end{theorem} 
Finally, in \cite{Togni}, Khennoufa and Togni determined the exact radio antipodal number for paths, resolving a long-standing conjecture.   
\begin{theorem}\cite{Togni} For $n \geq 5$,
\begin{equation*}
\an(P_n) = \left\{
\begin{array}{ll}
2m^2-2m+3 & \mbox{if $n = 2m+1$}, \\ [0.2cm]
2m^2-4m+5 & \mbox{if $n = 2m$}.
\end{array}
\right.
\end{equation*}
\end{theorem}
Observe that for $n = 2m+1$, $\binom{n-1}{2}-\frac{n-1}{2}+4$ is $2m^2-2m+4$, thus the upper bound in Theorem \ref{Pn:upper} is one more than the exact radio antipodal number of path $P_n$. Chartrand \emph{et al.} discussed the radio antipodal labeling of cycles, and determined $\an(C_n)$ for $n \equiv 2$ (mod $4$) and gave bounds for other cases in \cite{Zhang1}. In \cite{Juan}, Juan and Liu continued this study and determined $\an(C_n)$ for $n \equiv 1,3$ (mod $4$). They improved lower bound and gave an upper bound for $\an(C_n)$ for $n \equiv 0$ (mod $4$). Kumar and Panigrahi studied the optimal radio antipodal labeling of the Cartesian product of complete graphs and cycles in \cite{Kumar}. In \cite{Khennoufa}, Khennoufa and Togni studied the radio antipodal and radio number of the hypercube. William \emph{et al.} studied the radio antipodal numbers for graphs of  diameters $3$ and $4$ in \cite{William}.

In this paper, we determine the exact radio antipodal number for the Cartesian product of a complete graph and a path. The paper is organized as follows: In Section \ref{prel}, we give all necessary terms and definitions used in the present work. We also recall some existing results which we use to prove our main result. In Section \ref{main}, we first give an upper bound (Theorem \ref{thm:upper}) and next we provide a lower bound (Theorem \ref{thm:lower}) for the radio antipodal number of the Cartesian product of a complete graph and a path. Finally, we determine the exact radio antipodal number (Theorem \ref{thm:exact}) for the Cartesian product of a complete graph and a path. We illustrate our results with two examples (Examples \ref{exa1} and \ref{exa2}).  

\section{Preliminaries}
\label{prel}

We follow \cite{West} for standard graph-theoretic terms and notations. A complete graph $K_m$ is a graph with vertex set $V(K_m) = \{x_1,x_2,\ldots,x_m\}$, in which every pair of distinct vertices is adjacent. A path $P_n$ is a graph with vertex set $V(P_n) = \{y_1,y_2,\ldots,y_n\}$ and edge set $E(P_n) = \{y_iy_{i+1} : 1 \leq i \leq n-1\}$. The distance between two vertices $x$ and $y$ in $G$, denoted by $d_G(x,y)$, is the length of the shortest path joining $x$ and $y$ in $G$. The diameter of a graph $G$, denoted by $\diam(G)$ or simply $d$, is $\diam(G) = \max\{d_G(x,y) : x, y \in V(G)\}$. We drop the subscript in the above-defined terms when $G$ is clear in the context.   

Let $G=(V(G),E(G))$ and $H=(V(H),E(H))$ be two graphs. The Cartesian product of $G$ and $H$ is the graph $G \Box H$ with $V(G \Box H) = V(G) \times V(H)$ and two vertices $(a,b)$ and $(c,d)$ are adjacent if $a=c$ and $(b,d) \in E(H)$ or $b=d$ and $(a,c) \in E(G)$. 

\begin{observation}\label{Obs1} For a complete graph $K_m$ and a path $P_n$,
\begin{enumerate}[\rm (a)]
\item $|V(K_m \Box P_n)| = mn$,
\item $\diam(K_m \Box P_n) = n$,
\item $d_{K_m \Box P_n}((a,b),(c,d)) = d_{K_m}(a,c) + d_{P_n}(b,d)$.
\end{enumerate}
\end{observation}

Let $\vec{V} = (z_1,z_2,\ldots,z_{mn})$ be an ordering of $V(K_m \Box P_n)$. Then note that each $z_t\; (1 \leq t \leq mn)$ is an ordered pair $(x_{i_t},y_{j_t})$, where $x_{i_t} \in V(K_m)$ and $y_{j_t} \in V(P_n)$. Hence, each vertex of $K_m$ appears $n$ times and each vertex of $P_n$ appears $m$ times in $K_m \Box P_n$. For $z_a = (x_{i_a},y_{j_a}), z_b = (x_{i_b},y_{j_b}) \in V(K_m \Box P_n)\;(1 \leq a,b \leq mn)$, $x_{i_a}$ and $x_{i_b}$ are called \emph{distinct} if $i_a \neq i_b$. For any $z_a = (x_{i_a},y_{j_a}), z_b = (x_{i_b},y_{j_b}) \in V(K_m \Box P_n)\;(1 \leq a,b \leq mn)$, the distance between $z_a$ and $z_b$ is given by $d(z_a,z_b) = d(x_{i_a},x_{i_b})+d(y_{j_a},y_{j_b})$.

We recall the following two lemmas from \cite{Togni} and present them here with necessary adjustments in notation to align with the current work, which we use to derive the lower bound for the radio antipodal number of the Cartesian product of a complete graph and a path. 

\begin{lemma}\cite{Togni}\label{lemma1} Let $G$ be a graph of order $n$ and let $k$ be an integer. If $f$ is a radio $k$-labeling of $G$, then for any integer $k'> k$, there exists a radio $k'$-labeling $f'$ of $G$ with $\span(f') \leq \span(f) +(n-1)(k'-k)$. 
\end{lemma}

\begin{lemma} \cite{Togni}\label{lemma2} Let $G$ be a graph of order $n$ and let $k,k'$ be integers, $k'>k$. Given a radio $k$-labeling $f$ of $G$, let $z_1,z_2,\ldots,z_n$ be an ordering of the vertices of $G$ such that $f(z_i)\leq f(z_{i+1}),$ $1\leq i\leq n-1$ and let $\epsilon_i=|f(z_i)-f(z_{i-1})|-(1+k-d(z_i,z_{i-1})),$ $2\leq i\leq n$. Consider a set $I=\{i_1,i_2,\ldots,i_s\}\subset \{2,\ldots,n\},$ where $1\leq s\leq n-1$, such that $i_{j+1}>i_j+1$ for all $j,$ $1\leq j\leq s-1$. Then there exists a radio $k'$-labeling $f'$ of $G$ with $\span(f')\leq \span(f)+(n-1)(k'-k)-\sum_{i\in I} \min (k'-k,\epsilon_i)$. 
\end{lemma}

Using the above two lemmas, we establish the following relations between the radio antipodal number and the radio number of the graph $G$.

\begin{lemma}\label{lemma11} Let $G$ be a graph of order $n$ and diameter $d$. Then $\rn(G) \leq \an(G)+(n-1)$. 
\end{lemma}
\begin{proof} The proof follows by taking $k=d-1, k'=d$, and an optimal radio antipodal labeling $f$ in Lemma \ref{lemma1}. Since $\span(f) = \an(G)$ and $f'$ is a radio labeling of $G$, we have $\rn(G) \leq \span(f')$. Hence, by Lemma \ref{lemma1}, $\rn(G) \leq \span(f') \leq \span(f)+(n-1)(k'-k) = \an(G)+(n-1).$ \qed
\end{proof}

\begin{lemma}\label{lemma12} Let $G$ be a graph of order $n$ and diameter $d$. Given any optimal radio antipodal labeling $f$ of $G$, let $z_1,z_2,\ldots,z_n$ be an ordering of the vertices of $G$ such that $f(z_i)\leq f(z_{i+1}),$ $1 \leq i \leq n-1$ and let $\epsilon_i=|f(z_i)-f(z_{i-1})|-(d-d(z_i,z_{i-1})),$ $2\leq i\leq n$. Consider a set $I=\{i_1,i_2,\ldots,i_s\}\subset \{2,\ldots,n\},$ where $1 \leq s \leq n-1$, such that $i_{j+1} > i_j+1$ for all $j$, $1 \leq j \leq s-1$. Then $\rn(G) \leq \an(G)+(n-1)-\sum_{i \in I}\min(1,\epsilon_i)$.
\end{lemma}
\begin{proof} Again the proof follows by taking $k=d-1, k'=d$, and an optimal radio antipodal labeling $f$ in Lemma \ref{lemma2}. Since $\span(f) = \an(G)$ and $f'$ is a radio labeling of $G$, we have $\rn(G) \leq \span(f')$. Hence, by Lemma \ref{lemma2}, $\rn(G) \leq \span(f') \leq \span(f)+(n-1)(k'-k)-\sum_{i\in I} \min (k'-k,\epsilon_i) = \an(G)+(n-1)-\sum_{i \in I}\min(1,\epsilon_i)$. \qed
\end{proof}
   
\begin{lemma}\label{lemma3} Let $G$ be a graph of order $n$ and diameter $d$. Let $f$ be an optimal radio labeling defined on an ordering $z_1,z_2,\ldots,z_n$ of $V(G)$ by $f(z_1)=0$ and $f(z_{i+1})=f(z_i)+d+1-d(z_i,z_{i+1})$ for all $1 \leq i \leq n-1$. Then $\an(G) = \rn(G)-(n-1)$ if and only if there exists an optimal radio antipodal labeling $f'$ defined on the same ordering $z_1,z_2,\ldots,z_n$ of $V(G)$ such that $f'(z_i)=f(z_i)-(i-1)$ for all $1 \leq i \leq n$. 
\end{lemma}
\begin{proof} Suppose $\an(G)=\rn(G)-(n-1)$. Let $f$ be an optimal radio labeling described in the hypothesis, and let $f'$ be an optimal radio antipodal labeling defined on the same ordering. Since $\span(f') = \an(G) = \rn(G)-(n-1) = \span(f)-(n-1)$, we may assume that $f(z_1) = f'(z_1) = 0$. Thus, $f'(z_n) = f(z_n)-(n-1)$. Suppose, to the contrary, that $f'(z_i) \ne f(z_i)-(i-1)$ for some $i \in \{1,2,\ldots,n\}$. If $f'(z_a) < f(z_a)-(a-1)$ for some $a$ and $f'(z_i) \leq f(z_i)-(i-1)$ for all $i \neq a$. Then $f'(z_n) < f(z_n)-(n-1)$, which contradicts $\span(f') = \span(f)-(n-1)$. Similarly, if $f'(z_a) > f(z_a)-(a-1)$ for some $a$ and $f'(z_i) \geq f(z_i)-(i-1)$ for all $i \neq a$. Then $f'(z_n) > f(z_n)-(n-1)$, which also contradicts $\span(f') = \span(f)-(n-1)$. If there exist indices $a$ and $b$ such that $f'(z_a)>f(z_a)-(a-1)$ and $f'(z_b) < f(z_b)-(b-1)$. Then 
\begin{eqnarray*}
  f'(z_b)-f'(z_{b-1}) & \leq & f(z_b)-b-(f(z_{b-1})-(b-2)) \\
   &=& f(z_b)-f(z_{b-1})-2 \\
   &=&d+1-d(z_b,z_{b-1})-2 \\
   &=&d-d(z_b,z_{b-1})-1 \\
   &<&d-d(z_b,z_{b-1}),
\end{eqnarray*}
which contradicts the fact that $f'$ is a radio antipodal labeling of $G$. Hence, $f'(z_i)=f(z_i)-(i-1)$ for all $i\in \{1,2,\ldots,n\}$.

Conversely, suppose that an optimal radio antipodal labeling $f'$ defined on the ordering $z_1,z_2,\ldots,z_n$ of $V(G)$ satisfies  $f'(z_i)=f(z_i)-(i-1)$ for all $i \in \{1,2,\ldots,n\}$. Then, clearly $\an(G)=\span(f') = f'(z_n)=f(z_n)-(n-1)=\rn(G)-(n-1)$. \qed
\end{proof}

We use the following result of Kim \emph{et al.} \cite{Kim} on the radio number of the Cartesian product of a complete graph and a path to derive a lower bound for the radio antipodal number of this product. 

\begin{theorem}\cite{Kim}\label{Kim} Let $G=K_m \Box P_n$ be the Cartesian product of the complete graph $K_m$ and the path $P_n$ with $m\geq 3$ and $n\geq 4$. Then
\begin{equation*}\label{rn:exact}
\rn(G) = \left\{
\begin{array}{ll}
\frac{mn^2-2n+2}{2} & \mbox{if $n$ is even}, \\ [0.2cm]
\frac{mn^2-2n+m+2}{2} & \mbox{if $n$ is odd}.
\end{array}
\right.
\end{equation*}
\end{theorem}

\section{Main results}\label{main}
In this section, we continue to use terms and notations defined in the previous section. We first give an upper bound for the radio number of the Cartesian product of a complete graph and a path. 

\begin{theorem}\label{thm:upper} Let $m \geq 3$ and $n \geq 4$ be two integers. Then
\begin{equation*}\label{an:upper}
\an(K_m\Box P_n) \leq \left\{
\begin{array}{ll}
\frac{m(n-1)^2-2n+6}{2} & \mbox{if $n$ is odd}, \\ [0.2cm]
\frac{m(n^2-2n+4)-2n+4}{2} & \mbox{if $n$ is even}.
\end{array} 
\right.
\end{equation*}
\end{theorem}

\begin{proof} Denote the vertices of $K_m$ by $V(K_m)=\{x_1,x_2,\ldots,x_m\}$ and the vertices of $P_n$ by $V(P_n)=\{y'_1,y'_2,\ldots,y'_n\}$ with $E(P_n)=\{y'_iy'_{i+1} : i=1,2,\ldots,n-1\}$. We now consider the following two cases to define an ordering $z_1,z_2,\ldots,z_{mn}$ of $V(K
_m \Box P_n)$.

\textsf{Case-1:} $n$ is odd.

Define an ordering $y_1,y_2,\ldots,y_{n}$ of $V(P_n)$ as follows: 

If $n=5$, then set $y_1 := y'_4, y_2 := y'_1, y_3 := y'_3, y_4 := y'_5, y_5 := y'_2$.

If $n>5$, then set $y_1 := y'_{(n+3)/2}, y_2 := y'_2, y_3 := y'_n, y_4 := y'_{(n+1)/2}, y_5 := y'_1$ and for $6 \leq t \leq n$, set 
\begin{equation*}
y_t := 
\begin{cases}
y'_{(t+n-1)/2} & \mbox{if $t \equiv 0$ (mod $2$)}, \\
y'_{(t-1)/2} & \mbox{if $t \equiv 1$ (mod $2$)}. 
\end{cases}
\end{equation*}
From the above-defined ordering, if $n=5$, then $d(y_1,y_2)=d(y_4,y_5)=3$ and $d(y_2,y_3)=d(y_3,y_4)=d(y_1,y_5)=2$. Hence, $\sum_{j=1}^{4}d(y_j,y_{j+1})=10=(n^2-5)/2$. Moreover, $d(y_1,y_3)=d(y_3,y_5)=d(y_{4},y_1)=d(y_5,y_2)=1$ and $d(y_2,y_4)=4$.

If $n>5$, then
\begin{equation}\label{distancejtoj+1}
d(y_j,y_{j+1}) = 
\begin{cases}
\frac{n-1}{2} & \mbox{if $j=1,3,4$}, \\
n-2 & \mbox{if $j=2$},\\
\frac{n+3}{2} & \mbox{if $j=5$}, \\
\frac{n-1}{2} & \mbox{if $j\geq 6$ with $j \equiv 0$ (mod $2$)}, \\
\frac{n+1}{2} & \mbox{if $j\geq 7$ with $j \equiv 1$ (mod $2$)}, 
\end{cases}
\end{equation}
and 
\begin{equation}\label{distance1ton}
d(y_1,y_n) = 2.
\end{equation}
By \eqref{distancejtoj+1}, we get
\begin{eqnarray}\label{sumj=1ton-1}
 \sum_{j=1}^{n-1} d(y_j,y_{j+1}) & = &  \sum_{j=1}^{5} d(y_j,y_{j+1}) + \sum_{j=6}^{n-1} d(y_j,y_{j+1}) \non \\
 & = &
 3\(\frac{n-1}{2}\)+n-2+\frac{n+3}{2}+\(\frac{n-5}{2}\)\(\frac{n-1}{2}\)+\(\frac{n-7}{2}\)\(\frac{n+1}{2}\) \non \\
 & = & \frac{n^2-5}{2}.
\end{eqnarray}

Moreover,
\begin{equation}\label{distancejtoj+2}
d(y_j,y_{j+2}) = 
\begin{cases}
\frac{n-3}{2} & \mbox{if $j=1,2$}, \\
n-1 & \mbox{if $j=3$},\\
2 & \mbox{if $j=4,5$}, \\
1 & \mbox{if $6\leq j\leq n-2$}, 
\end{cases}
\end{equation}

and

\begin{equation}\label{distance1ton-1}
d(y_{n-1},y_1) = d(y_n,y_2) = \frac{n-5}{2}.
\end{equation}

We now consider the following two sub-cases in order to define an ordering $z_1,z_2,\ldots,z_{mn}$ of $V(K_m \Box P_n)$.

\textsf{Subcase-1:} $m\geq n$.

For $1 \leq t \leq mn$, define $z_t := (x_i,y_j)$, where $i \in \{1,2,\ldots,m\}$, $j \in \{1,2,\ldots,n\}$, and  

$$t := \left\{
\begin{array}{ll}
n(i-j)+j & \mbox{if $i\geq j$}, \\[0.2cm]
n(m+i-j)+j & \mbox{if $i<j$}.
\end{array}
\right.$$ 

\textsf{Subcase-2:} $m<n$.

For $1 \leq t \leq mn$, define an ordering $z_1,z_2,\ldots,z_{mn}$ of $V(K_m \Box P_n)$ as follows: 

For $1\leq t\leq \lcm(m,n)$, set 
$$z_t := \left\{
\begin{array}{ll}
(x_i,y_j) & \mbox{if $t \equiv i$ (mod $m$), $t \equiv j$ (mod $n$) and $1\leq i\leq m-1$, $1\leq j\leq n-1$}, \\[0.2cm]
(x_i,y_n) & \mbox{if $t \equiv i$ (mod $m$), $t \equiv j$ (mod $n$) and $1\leq i\leq m-1$, $j=0$},  \\[0.2cm]
(x_m,y_j) & \mbox{if $t \equiv i$ (mod $m$), $t \equiv j$ (mod $n$) and $i=0$, $1\leq j\leq n-1$}, \\[0.2cm]
(x_m,y_n) & \mbox{if $t \equiv i$ (mod $m$), $t \equiv j$ (mod $n$) and $i=0$, $j=0$}. \\[0.2cm]
\end{array}
\right.$$ 

For $\lcm (m,n)+1\leq t\leq 2 \cdot \lcm (m,n)$, set 
$$z_t := \left\{
\begin{array}{ll}
(x_{i+1},y_j) & \mbox{  if  $t \equiv i$ (mod $m$), $t \equiv j$ (mod $n$) and $1\leq j\leq n-1$}, \\[0.2cm]
(x_{i+1},y_n) & \mbox{  if  $t \equiv i$ (mod $m$), $t \equiv j$ (mod $n$) and $j=0$}. \\[0.2cm]
\end{array}
\right.$$

For $r\cdot\lcm(m,n)+1 \leq t \leq (r+1)\cdot\lcm (m,n)$, where $2 \leq r \leq \gcd(m,n)-1$, set
\begin{equation*}
z_t := \left\{
\begin{array}{ll}
(x_{i+r},y_j) & \mbox{if $t \equiv i$ (mod $m$), $t \equiv j$ (mod $n$) and $0\leq i\leq m-r$, $1\leq j\leq n-1$}, \\[0.2cm]
(x_{i+r},y_n) & \mbox{if $t \equiv i$ (mod $m$), $t \equiv j$ (mod $n$) and $0\leq i\leq m-r$, $j=0$}, \\[0.2cm]
(x_{i+r-m},y_j) & \mbox{if $t \equiv i$ (mod $m$), $t \equiv j$ (mod $n$) and $m-r+1\leq i\leq m-1$, $1\leq j\leq n-1$}, \\[0.2cm]
(x_{i+r-m},y_n) & \mbox{if $t \equiv i$ (mod $m$), $t \equiv j$ (mod $n$) and $m-r+1\leq i\leq m-1$, $j=0$}. 
\end{array}
\right.
\end{equation*}

For the above-defined ordering, it is clear that $x_{i_t}$, $x_{i_{t+1}}$, and $x_{i_{t+2}}$ are distinct for $1 \leq t \leq mn-2$. Hence, $d(x_{i_t}, x_{i_{t+1}})=d(x_{i_t}, x_{i_{t+2}})=1$ for all $1 \leq t \leq mn-2$. By \eqref{distancejtoj+1}, \eqref{distance1ton}, \eqref{distancejtoj+2} and \eqref{distance1ton-1} we get
\begin{equation}\label{gap2}
\sum_{t=a}^{a+1} d(y_{j_t}, y_{j_{t+1}}) \leq n-1+d(y_{j_a}, y_{j_{a+2}}),
\end{equation}
and 
\begin{equation}\label{gap3}
\sum_{t=a}^{a+2} d(y_{j_t}, y_{j_{t+1}}) \leq 
\begin{cases}
8 & \mbox{  if  $n=5$}, \\
2n-3 & \mbox{  if  $n>5$}.
\end{cases}
\end{equation}

By \eqref{distance1ton} and \eqref{sumj=1ton-1} we get
\begin{eqnarray}\label{yjtsum:odd}
\sum_{t=1}^{mn-1} d(y_{j_t},y_{j_{t+1}}) & = & m \sum_{j=1}^{n-1} d(y_j,y_{j+1}) + (m-1) d(y_1,y_n) \non \\
& = & \frac{m(n^2-5)}{2}+2(m-1) \non \\ 
& = & \frac{mn^2-m-4}{2}.
\end{eqnarray}

Define a function $f$ on an ordering $z_1,z_2,\ldots,z_{mn}$ of $V(K_m \Box P_n)$ by $f(z_1)=0$, and for $1 \leq t \leq mn-1$, $f(z_{t+1})=f(z_t)+n-d(z_t,z_{t+1})$. 

\textsf{Claim-1:} $f$ is a radio antipodal labeling of $K_m\Box P_n$.

Let $z_a$ and $z_b$ be two arbitrary vertices of $V(K_m\Box P_n)$ with $1\leq a<b \leq mn$. Then
\begin{eqnarray}\label{zazb}
f(z_b)-f(z_a) & = &  \sum_{t=a}^{b-1} (f(z_{t+1})-f(z_t)) \non \\
& = & \sum_{t=a}^{b-1} (n-d(z_t,z_{t+1})) \non \\
& = & \sum_{t=a}^{b-1} (n-d((x_{i_t},y_{j_t}),(x_{i_{t+1}},y_{j_{t+1}}))) \non \\
& = & (b-a)n-\sum_{t=a}^{b-1} (d(x_{i_t},x_{i_{t+1}})+d(y_{j_t},y_{j_{t+1}})) \non \\
& = & (b-a)n-(b-a)-\sum_{t=a}^{b-1} d(y_{j_t},y_{j_{t+1}}) \non \\
& = & (b-a)(n-1)-\sum_{t=a}^{b-1} d(y_{j_t},y_{j_{t+1}}).
\end{eqnarray}

If $b-a=1$, then by the above-defined labeling, $f(z_b)-f(z_a)=n-d(z_a,z_b)$, $1\leq a<b \leq mn$. 

If $b-a=2$, then by \eqref{gap2} and \eqref{zazb}, $f(z_b)-f(z_a)=2(n-1)-\sum_{t=a}^{a+1} d(y_{j_t},y_{j_{t+1}})\geq 2(n-1)+1-n-d(y_{j_a},y_{j_{a+2}})=n-1-d(y_{j_a},y_{j_{a+2}})=n-d(x_{i_a},x_{i_{a+2}})-d(y_{j_a},y_{j_{a+2}})=n-d((x_{i_a},y_{j_a}),(x_{i_{a+2}},y_{j_{a+2}}))=n-d(z_a,z_b)$.

If $b-a=3$, then by \eqref{gap3} and \eqref{zazb}, $f(z_b)-f(z_a)=3(n-1)-\sum_{t=a}^{a+2} d(y_{j_t},y_{j_{t+1}})$. If $n=5$, then $f(z_b)-f(z_a)\geq 12-8=4=n-1\geq n-d(z_a,z_b)$. If $n>5$, then $f(z_b)-f(z_a) \geq 3(n-1)-(2n-3)=n\geq n-d(z_a,z_b)$. 

If $b-a>3$, then $f(z_b)-f(z_a)=f(z_b)-f(z_{a+3})+f(z_{a+3})-f(z_a)\geq n-1+f(z_b)-f(z_{a+3})\geq n>n-d(z_a,z_b)$, since  $d(z_a,z_b) \geq 1$. Therefore, in all cases, $f$ is a radio antipodal labeling of $K_m\Box P_n$. This completes the proof of Claim-1.

The span of $f$ is
\begin{eqnarray*}
    \span(f) & = & f(z_{mn})-f(z_1) \\
    & = & \sum_{t=1}^{mn-1} (f(z_{t+1})-f(z_t)) \\
    & = & \sum_{t=1}^{mn-1} (n-d(z_t,z_{t+1}))\\
    & = & (mn-1)n-\sum_{t=1}^{mn-1} (d(x_{i_t},x_{i_{t+1}})+d(y_{j_t},y_{j_{t+1}})) \\
    & = & (mn-1)(n-1)-\sum_{t=1}^{mn-1} d(y_{j_t},y_{j_{t+1}}) \\
    & = & (mn-1)(n-1)-\(\frac{mn^2-m-4}{2}\) \hspace{1cm}\text{(using \eqref{yjtsum:odd})} \\
    & = & \frac{m(n-1)^2-2n+6}{2}.
\end{eqnarray*}

\textsf{Case-2:} $n$ is even.

Define an ordering $y_1,y_2,\ldots,y_{n}$ as follows: Set $y_1 := y'_{(n+2)/2}, y_2 := y'_1, y_{n-1} := y'_n, y_n := y'_{n/2}$ and for $3 \leq t \leq n-2$, set 
\begin{equation*}
y_t := 
\begin{cases}
y'_{(n-t+2)/2} & \mbox{if $t \equiv 0$ (mod $2$)}, \\
y'_{(2n-t+1)/2} & \mbox{if $t \equiv 1$ (mod $2$)}.
\end{cases}
\end{equation*}

From the above-defined ordering, we have  

\begin{equation}\label{distancejtoj+1even}
d(y_j,y_{j+1}) = 
\begin{cases}
\frac{n}{2} & \mbox{if $j \equiv 1$ (mod $2$)}, \\
n-2 & \mbox{if $j=2$ or $n-2$},\\
\frac{n-2}{2} & \mbox{if $j \equiv 0$ (mod $2$) and $4\leq j \leq n-4$},
\end{cases}
\end{equation}

and 

\begin{equation}\label{distance1toneven}
d(y_1,y_n) = 1.
\end{equation}

By \eqref{distancejtoj+1even}, we get
\begin{eqnarray}\label{sumj=1ton-1even}
\sum_{j=1}^{n-1} d(y_j,y_{j+1}) & = & \(\frac{n}{2}\)\(\frac{n}{2}\)+2(n-2)+\(\frac{n-6}{2}\)\(\frac{n-2}{2}\) \non \\
& = & \frac{n^2-2}{2}.
\end{eqnarray}

Moreover,
\begin{equation}\label{distancejtoj+2even}
d(y_j,y_{j+2}) = 
\begin{cases}
\frac{n-4}{2} & \mbox{if $j=1,2,n-3,n-2$ }, \\
1 & \mbox{if $3\leq j\leq n-4$}, 
\end{cases}
\end{equation}

and

\begin{equation}\label{distance1ton-1even}
d(y_{n-1},y_1) = d(y_n,y_2) = \frac{n-2}{2}.
\end{equation}

Define an ordering $z_1,z_2,\ldots,z_{mn}$ of the vertex set  $V(K_m \Box P_n)$ as described in Case 1, depending on whether $m\geq n$ and $m<n$.

Then it is clear that $x_{i_t}$, $x_{i_{t+1}}$ and $x_{i_{t+2}}$ are distinct for $1 \leq t \leq mn-2$. Hence, $d(x_{i_t}, x_{i_{t+1}})=d(x_{i_t}, x_{i_{t+2}})=1$ for all $1 \leq t \leq mn-2$. By \eqref{distancejtoj+1even}, \eqref{distance1toneven}, \eqref{distancejtoj+2even} and \eqref{distance1ton-1even}, we obtain 

\begin{equation}\label{gap2even}
\sum_{t=a}^{a+1} d(y_{j_t}, y_{j_{t+1}}) = 
\begin{cases}
n+d(y_{j_a}, y_{j_{a+2}}) & \mbox{if  $a\equiv 1,2,n-3,n-2$ (mod $n$)}, \\
\frac{n+2}{2} & \mbox{if  $a\equiv n-1$ or $0$ (mod $n$)} \\
n-1 & \mbox{otherwise},
\end{cases}
\end{equation}
and 
\begin{equation}\label{gap3even}
\sum_{t=a}^{a+2} d(y_{j_t}, y_{j_{t+1}}) \leq 
\begin{cases}
5 & \mbox{if $n=4$ and $a \equiv 3$ (mod $4$)}, \\
7 & \mbox{if $n=4$ and $a \equiv 0,1,2$ (mod $4$)}, \\
11 & \mbox{if  $n=6$}, \\
2n-2 & \mbox{if  $n>6$}.
\end{cases}
\end{equation}

By \eqref{distance1toneven} and \eqref{sumj=1ton-1even}, we get
\begin{eqnarray}\label{yjtsum:even}
\sum_{t=1}^{mn-1} d(y_{j_t},y_{j_{t+1}}) & = & m \sum_{j=1}^{n-1} d(y_j,y_{j+1}) + (m-1) d(y_1,y_n) \non \\
& = & \frac{m(n^2-2)}{2}+(m-1) \non \\ 
& = & \frac{mn^2-2}{2}.
\end{eqnarray}

Define a function $f$ on an ordering $z_1,z_2,\ldots,z_t$ of $V(K_m \Box P_n)$ by $f(z_1)=0$ and for $1 \leq t\leq mn-1$, 

$$f(z_{t+1}) = \left\{
\begin{array}{ll}
f(z_t)+n+1-d(z_t,z_{t+1}) & \mbox{if  $t\equiv 2$ or $n-2$ (mod $n$)}, \\[0.2cm]
f(z_t)+n-d(z_t,z_{t+1}) & \mbox{otherwise}. \\[0.2cm]
\end{array}
\right.
$$

\textsf{Claim-2:} $f$ is a radio antipodal labeling of $K_m\Box P_n$.

Let $z_a$ and $z_b$ be two arbitrary vertices of $V(K_m\Box P_n)$ with $1\leq a<b \leq mn$. Then
\begin{equation}\label{zazbeven}
f(z_b)-f(z_a)\geq (b-a)(n-1)-\sum_{t=a}^{b-1} d(y_{j_t},y_{j_{t+1}}).
\end{equation}

If $b-a=1$, then by above-defined labeling, $f(z_b)-f(z_a)\geq n-d(z_a,z_b)$ for all $1\leq a<b \leq mn$.

If $b-a=2$, then for $a \equiv 1,2,n-3,n-2$ (mod $n$), $f(z_b)-f(z_a)=f(z_{a+2})-f(z_{a+1})+f(z_{a+1})-f(z_a)=n+1-d(z_a,z_{a+1})+n-d(z_a,z_{a+1})=2n+1-\sum_{t=a}^{a+1} d(z_{t}, z_{t+1})=2n+1-\sum_{t=a}^{b-1} (d(x_{i_t},x_{i_{t+1}})+d(y_{j_t},y_{j_{t+1}}))=2n+1-2-\sum_{t=a}^{a+1} d(y_{j_t}, y_{j_{t+1}})$. Then by \eqref{gap2even}, $f(z_b)-f(z_a)=2n-1-n-d(y_{j_a}, y_{j_{a+2}})=n-d(x_{i_a}, x_{i_{a+2}})-d(y_{j_a}, y_{j_{a+2}})=n-d(z_a,z_b)$. For $a\equiv n-1$ or $0$ (mod $n$), $f(z_b)-f(z_a)=2(n-1)-\sum_{t=a}^{a+1} d(y_{j_t},y_{j_{t+1}})= 2(n-1)-(n+2)/2=3(n-2)/2 \geq n-1 \geq n-d(z_a,z_b)$ as $d(z_a,z_b)\geq 2$. For remaining cases, $f(z_b)-f(z_a)=2(n-1)-\sum_{t=a}^{a+1} d(y_{j_t},y_{j_{t+1}})= 2(n-1)-(n-1)=n-1\geq n-d(z_a,z_b)$ as $d(z_a,z_b)>1$. 

Assume $b-a = 3$. If $n = 4$, then for $a\equiv 0,1$ or $2$ (mod $4$), $f(z_b)-f(z_a)=5+4+4-\sum_{t=a}^{a+2} d(z_{t}, z_{{t+1}})=13-\sum_{t=a}^{a+2} (d(x_{i_t},x_{i_{t+1}})+d(y_{j_t},y_{j_{t+1}}))=13-3-\sum_{t=a}^{a+2} d(y_{j_t}, y_{j_{t+1}})$. Then by \eqref{gap3even}, $f(z_b)-f(z_a)=10-\sum_{t=a}^{a+2} d(y_{j_t}, y_{j_{t+1}})\geq 10-7=3\geq 4-d(z_a,z_b)$ as $d(z_a,z_b)\geq 1$. For $a\equiv 3$ (mod $n$), $f(z_b)-f(z_a)=3(4)-\sum_{t=a}^{a+2} d(z_{t},z_{t+1})=12-3-\sum_{t=a}^{a+2} d(y_{j_t}, y_{j_{t+1}})$. Then by \eqref{gap3even}, $f(z_b)-f(z_a)= 9-5=4\geq 4-d(z_a,z_b)$. If $n=6$, then by \eqref{gap3even} and \eqref{zazbeven}, $f(z_b)-f(z_a)\geq 3(n-1)-\sum_{t=a}^{a+2} d(y_{j_t},y_{j_{t+1}}) \geq 3n-3-11=4=n-2 \geq n-d(z_a,z_b)$ as $d(z_a,z_b)\geq 2$. If $n > 6$, then by \eqref{gap3even} and \eqref{zazbeven}, $f(z_b)-f(z_a)\geq 3(n-1)-\sum_{t=a}^{a+2} d(y_{j_t},y_{j_{t+1}}) \geq 3(n-1)-(2n-2)=n-1\geq n-d(z_a,z_b)$.  

If $b-a>3$, then $f(z_b)-f(z_a)=f(z_b)-f(z_{a+3})+f(z_{a+3})-f(z_a)\geq n-2+f(z_b)-f(z_{a+3})\geq n-1\geq n-d(z_a,z_b)$ as $d(z_a,z_b)\geq 1$. Therefore, in all cases, $f$ is a radio antipodal labeling of $K_m\Box P_n$. This completes the proof of Claim-2. 

If $n = 4$, then the span of $f$ is
\begin{eqnarray*}
    \span(f) & = & f(z_{4m})-f(z_1) \\
    & = & \sum_{t=1}^{4m-1} (f(z_{t+1})-f(z_t)) \\
    & = & \sum_{\substack{t=1 \\ t \not\equiv 2 (\text{mod } 4)}}^{4m-1}(f(z_{t+1})-f(z_t))+\sum_{\substack{t=1 \\ t \equiv 2 (\text{mod } 4)}}^{4m-1}(f(z_{t+1})-f(z_t)) \\
    & = & \sum_{\substack{t=1 \\ t \not\equiv 2 (\text{mod } 4)}}^{4m-1}(4-d(z_t,z_{t+1}))+\sum_{\substack{t=1 \\ t \equiv 2 (\text{mod } 4)}}^{4m-1}(5-d(z_t,z_{t+1})) \\
    & = & \sum_{t=1}^{4m-1} (4-d(z_t,z_{t+1}))+m\\
    & = & 4(4m-1)-\sum_{t=1}^{4m-1} (d(x_{i_t},x_{i_{t+1}})+d(y_{j_t},y_{j_{t+1}}))+m\\
    & = & 3(4m-1)-\sum_{t=1}^{4m-1} d(y_{j_t},y_{j_{t+1}})+m \\
    & = & 3(4m-1)-(8m-1)+m \hspace{1cm}\text{(using \eqref{yjtsum:even})} \\
    & = & 5m-2.
    \end{eqnarray*}

If $n \geq 6$, then the span of $f$ is 
\begin{eqnarray*}
    \span(f) & = & f(z_{mn})-f(z_1) \\
    & = & \sum_{t=1}^{mn-1} (f(z_{t+1})-f(z_t)) \\
    & = & \sum_{\substack{t=1 \\ t \not\equiv 2, n-2(\text{mod } n)}}^{4m-1}(f(z_{t+1})-f(z_t))+\sum_{\substack{t=1 \\ t \equiv 2, n-2(\text{mod } n)}}^{4m-1}(f(z_{t+1})-f(z_t)) \\
    & = & \sum_{\substack{t=1 \\ t \not\equiv 2, n-2(\text{mod } n)}}^{4m-1}(n-d(z_t,z_{t+1}))+\sum_{\substack{t=1 \\ t \equiv 2, n-2(\text{mod } n)}}^{4m-1}(n+1-d(z_t,z_{t+1})) \\
    & = & \sum_{t=1}^{mn-1} (n-d(z_t,z_{t+1}))+2m\\
    & = & (mn-1)n-\sum_{t=1}^{mn-1} (d(x_{i_t},x_{i_{t+1}})+d(y_{j_t},y_{j_{t+1}}))+2m\\
    & = & (mn-1)(n-1)-\sum_{t=1}^{mn-1} d(y_{j_t},y_{j_{t+1}})+2m \\
    & = & (mn-1)(n-1)-\(\frac{mn^2-2}{2}\)+2m \hspace{1cm}\text{(using \eqref{yjtsum:even})} \\
    & = & \frac{m(n^2-2n+4)-2n+4}{2},
\end{eqnarray*}
which completes the proof. \qed
\end{proof}

We now give a lower bound for the radio number of the Cartesian product of a complete graph and a path. 

\begin{theorem}\label{thm:lower} Let $m\geq 3$ and $n \geq 4$ be two integers. Then
\begin{equation*}\label{an:lower}
\an(K_m\Box P_n) \geq \left\{
\begin{array}{ll}
\frac{m(n-1)^2-2n+4}{2} & \mbox{if $n$ is odd}, \\ [0.2cm]
\frac{m(n^2-2n+4)-2n+4}{2} & \mbox{if $n$ is even}.
\end{array}
\right.
\end{equation*}
\end{theorem}

\begin{proof} We consider the following two cases.

\textsf{Case-1:} $n$ is odd.

By Lemma \ref{lemma11}, we have $\rn (K_m\Box P_n)\leq \an (K_m\Box P_n)+(mn-1)$, and by Theorem \ref{Kim}, we have $\rn (K_m\Box P_n)=\frac{mn^2-2n+m+2}{2}$. Using these, we get
\begin{eqnarray*}
\an (K_m\Box P_n) & \geq &  \rn (K_m\Box P_n)-(mn-1) \\
& = & \(\frac{mn^2-2n+m+2}{2}\)-(mn-1) \\
& = & \frac{m(n-1)^2-2n+4}{2}.
\end{eqnarray*}

\textsf{Case-2:} $n$ is even.

Recall that $d=\diam (K_m\Box P_n)=n$. We consider Lemma \ref{lemma12} with the radio antipodal labeling $f$ of $K_m\Box P_n$ described in the proof of Theorem \ref{thm:upper} with $k=d-1=n-1$ and $k'=d=n$. By Lemma \ref{lemma12}, we have $\epsilon_t=|f(z_t)-f(z_{t-1})|-(n-d(z_t,z_{t-1}))$. As per defined radio antipodal labeling of $z_t:=(x_{i_t},y_{j_t})$, $2 \leq t \leq mn$,
\begin{equation*}
\epsilon_t = \left\{
\begin{array}{ll}
1 & \mbox{if  $t\equiv 3$ or $n-1$ (mod $n$)}, \\ [0.2cm]
0 & \mbox{otherwise}.
\end{array}
\right.
\end{equation*}

Observe that each $y_j$ appears $m$ times in $K_m\Box P_n$. If $n = 4$, then for one value of $j$, the corresponding value of $\epsilon_t$ is 1. This gives $\sum_{t \in I}\min(1, \epsilon_t) = m$, since $\epsilon_t=1$ for all $t \in I=\{qn+3 : 0 \leq q \leq m-1\}$.

By Lemma \ref{lemma12}, we have $\rn (K_m\Box P_4)\leq \an (K_m\Box P_4)+(4m-1)-\sum_{t\in I} \min (1,\epsilon_t)$, and by Theorem \ref{Kim}, we have $\rn (K_m\Box P_4)=8m-3$. Using these, we obtain
\begin{eqnarray*}
\an (K_m\Box P_4) & \geq &  \rn (K_m\Box P_4)-(4m-1)+\sum_{t\in I} \min (1,\epsilon_t) \\
& = & 8m-3-(4m-1)+m \\
& = & 5m-2. 
\end{eqnarray*}

If $n > 4$, then for two values of $j$, the corresponding values of $\epsilon_t$ are $1$. This gives $\sum_{t\in I} \min (1,\epsilon_t)=2m$, since $\epsilon_t=1$ for all $t \in I$, where $I = \{qn+3,qn+n-1 : 0 \leq q \leq m-1\}$.

By Lemma \ref{lemma12}, we have $\rn (K_m\Box P_n)\leq \an (K_m\Box P_n)+(mn-1)-\sum_{t\in I} \min(1,\epsilon_t)$, and by Theorem \ref{Kim}, we have $\rn (K_m\Box P_n)=\frac{mn^2-2n+2}{2}$. Using these, we obtain
\begin{eqnarray*}
\an (K_m\Box P_n) & \geq &  \rn (K_m\Box P_n)-(mn-1)+\sum_{t\in I} \min(1,\epsilon_t) \\
& = & \(\frac{mn^2-2n+2}{2}\)-(mn-1)+2m \\
& = & \frac{m(n^2-2n+4)-2n+4}{2}. 
\end{eqnarray*} \qed
\end{proof}

We now determine the exact radio antipodal number of the Cartesian product of a complete graph and a path in the following result. 

\begin{theorem}\label{thm:exact} Let $m, n \geq 3$ be two integers. Then
\begin{equation*}\label{an:exact}
\an(K_m\Box P_n) = \left\{
\begin{array}{ll}
6 & \mbox{if $m = n = 3$}, \\ [0.2cm]
2m-1 & \mbox{if $m \geq 4$ and $n = 3$}, \\ [0.2cm]
5m-2 & \mbox{if $m \geq 3$ and $n = 4$}, \\ [0.2cm]
\frac{m(n-1)^2-2n+6}{2} & \mbox{if $m \geq 3$ and $n > 4$ is odd}, \\ [0.2cm]
\frac{m(n^2-2n+4)-2n+4}{2} & \mbox{if $m \geq 3$ and $n > 4$ is even}.
\end{array}
\right.
\end{equation*}
\end{theorem}
\begin{proof} We consider the following cases. 

\textsf{Case-1:} $m = n = 3$. 

It is easy to verify that $\an(K_3 \Box P_3) = 6$. An optimal radio antipodal labeling of $K_3 \Box P_3$ is shown in following Figure \ref{K3P3}.  
\begin{figure}[ht]
\centering
\includegraphics[width=0.5\linewidth]{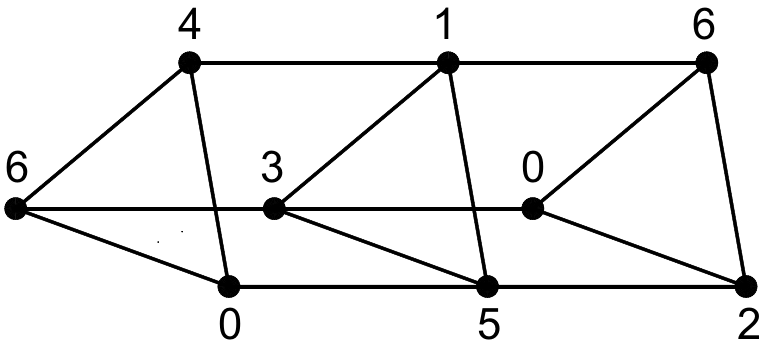}
\caption{An optimal radio antipodal labeling of $K_3 \Box P_3$ with $\an(K_3 \Box P_3) = 6$.}
\label{K3P3}
\end{figure}

\textsf{Case-2:} $m \geq 4$ and $n = 3$.

Denote $V(K_m) = \{x_1,x_2,\ldots,x_m\}$ and $V(P_3) = \{y_1,y_2,y_3\}$ and $E(P_3) = \{y_1y_2,y_2y_3\}$. Let $f$ be an optimal radio antipodal labeling of the graph $K_m \Box P_3$. Let $R_1 = \{(x_i,y_1) : 1 \leq i \leq m \}$, $R_2 = \{(x_i,y_2) : 1 \leq i \leq m \}$ and $R_3 = \{(x_i,y_3) : 1 \leq i \leq m\}$.  Note that 
$d((x_{i_1},y_1), (x_{i_2},y_1)) = d((x_{j_1},y_2), (x_{j_2},y_2)) = 1$ for all $1 \leq i_1, i_2, j_1,j_2 \leq m$ with $i_1 \neq i_2, j_1 \neq j_2$. It follows that 
$|f(R_1)| = |f(R_2)| = m$,
meaning that all vertices of $R_1$ and $R_2$ must receive distinct labels. Furthermore, since 
$d((x_{i},y_1), (x_{i},y_2)) = 1$ for each $1 \leq i \leq m$ and $d((x_i,y_1), (x_j,y_2)) = 2$
for all $1 \leq i, j \leq m$ with $i \neq j$, the sets $f(R_1)$ and $f(R_2)$ must be disjoint. Therefore, the total number of distinct labels assigned to $R_1 \cup R_2$ is $2m$, implying that the largest label used is at least $2m - 1$ as labeling starts with 0. Now, we define the function $f : \{(x_i,y_j) : 1\leq i\leq m, 1\leq j\leq 3\} \rightarrow \{0,1,2,\ldots,2m-1\}$ as follows:
\begin{equation*}
f((x_i,y_j)) = \left\{
\begin{array}{ll}
2(i-1) & \mbox{if $1\leq i\leq m$ and $j=1$}, \\ [0.2cm]
2(m+i-2)-1 & \mbox{if $1 \leq i \leq 2$ and $j = 2$}, \\ [0.2cm]
2(i-2)-1 & \mbox{if $3 \leq i \leq m$ and $j = 2$}, \\ [0.2cm]
2(m+i-4) & \mbox{if $1 \leq i \leq 3$ and $j = 3$}, \\ [0.2cm]
2(i-4) & \mbox{if $4 \leq i \leq m$ and $j = 3$}.
\end{array} 
\right.
\end{equation*}
We show that the function $f$ defined above is a radio antipodal labeling of $K_m \Box P_3$ for $m \geq 4$. Let $(x_a,y_b)$ and $(x_c,y_d)$ be two arbitrary vertices of $K_m \Box P_3$. If $(x_a,y_b),(x_c,y_d) \in R_1$, then note that $d((x_a,y_b),(x_c,y_d)) = 1$. Therefore 
\begin{eqnarray*}
  |f((x_a,y_b))-f((x_c,y_d))| &=& |2(a-1)-2(c-1)| \\
  & = & 2|a-c| \\
  & \geq & 2 = d-d((x_a,y_b),(x_c,y_d)).
\end{eqnarray*}
Hence, $f$ satisfies the radio antipodal labeling condition for any two vertices in $R_1$. Similarly, one can verify that $f$ satisfies the radio antipodal labeling condition for any two vertices in $R_2$ or in $R_3$. Now, suppose that $(x_a,y_b) \in R_1$ and $(x_c,y_d) \in R_2$, then note that $d((x_a,y_b),(x_c,y_d)) = 1$ if $a=c$, and $d((x_a,y_b),(x_c,y_d)) = 2$ if $a \neq c$. First, suppose that $a = c$ and $a \in \{1,2\}$. Then 
\begin{eqnarray*}
|f((x_a,y_b))-f((x_c,y_d))| & = & |2(a-1)-2(m+a-2)+1| \\
& = & |3-2m| \\
& \geq & 5 > d-d((x_a,y_b),(x_c,y_d)).
\end{eqnarray*}
Next, suppose that $a = c$ and $a \in \{3,4,\ldots,m\}$. Then 
\begin{eqnarray*}
|f((x_a,y_b))-f((x_c,y_d))| & = & |2(a-1)-2(a-2)+1| \\
& = & 3 > d-d((x_a,y_b),(x_c,y_d)).
\end{eqnarray*}
Suppose now that $a \neq c$ and $c \in \{1,2\}$. Then 
\begin{eqnarray*}
|f((x_a,y_b))-f((x_c,y_d))| & = & |2(a-1)-2(m+c-2)+1| \\
& = & |3+2(a-m-c)| \\
& \geq & 1 \geq d-d((x_a,y_b),(x_c,y_d)).
\end{eqnarray*}
Finally, suppose that $a \neq c$ and $c \in \{3,4,\ldots,m\}$. Then 
\begin{eqnarray*}
|f((x_a,y_b))-f((x_c,y_d))| & = & |2(a-1)-2(c-2)+1| \\
& = & |3+2(a-c)| \\
& \geq & 1 \geq d-d((x_a,y_b),(x_c,y_d)).
\end{eqnarray*}
Thus, $f$ satisfies the radio antipodal condition for every pair of vertices with one vertex in $R_1$ and the other in $R_2$. Similarly, one can verify the condition for pairs with one vertex in $R_1$ and the other in $R_3$, and for pairs with one vertex in $R_2$ and the other in $R_3$. Hence, $f$ is a radio antipodal labeling of $K_m \Box P_3$. 

\textsf{Case-3:} $m \geq 3$ and $n > 4$ is odd,

By Theorem \ref{thm:lower}, the lower bound for $\an(K_m\Box P_n)$ is $\frac{m(n-1)^2-2n+4}{2}$. Suppose that $\an(K_m\Box P_n)$ $= \frac{m(n-1)^2-2n+4}{2} = \rn(K_m\Box P_n)-(mn-1)$. Consider an optimal radio labeling $f$ of $K_m \Box P_n$ described in \cite{Kim}, and let $z_1,z_2,\ldots,z_{mn}$ be the ordering of $V(K_m \Box P_n)$ induced by $f$. Note that $f$ is defined on $z_1,z_2,\ldots,z_{mn}$ by $f(z_1)=0$ and $f(z_{i+1}) = f(z_i)+n+1-d(z_i,z_{i+1})$. Then by Lemma \ref{lemma3}, there exists an optimal radio antipodal labeling $f'$ defined on the ordering $z_1,z_2,\ldots,z_n$ of $V(G)$ such that $f'(z_i)=f(z_i)-(i-1)$ for all $i\in \{1,2,\ldots,n\}$. Now observe that for an ordering $z_1,z_2,\ldots,z_{mn}$ given in \cite{Kim}, there is a vertex $z_i$ with $i=n(m-1)$ such that $d(z_{i-1},z_i)=(n+3)/2$, $d(z_{i},z_{i+1})=n-1$ and $d(z_{i-1},z_{i+1})=(n-1)/2$.  Then we have
\begin{eqnarray*}
f'(z_{i+1})-f'(z_{i-1}) & = & f(z_{i+1})-i-f(z_{i-1})+i-2  \\
& = & f(z_{i+1})-f(z_i)+f(z_i)-f(z_{i-1})-2 \\
& = & 2(n+1)-\(\frac{n+3}{2}\)-(n-1)-2 \\
& = & \frac{n-1}{2}.
\end{eqnarray*}
Note that $n-d(z_{i-1},z_{i+1})=n-(n-1)/2=(n+1)/2>(n-1)/2=f'(z_{i+1})-f'(z_{i-1})$, which is a contradiction as $f'$ is a radio antipodal labeling. Therefore, our assumption $\an(K_m\Box P_n) = \frac{m(n-1)^2-2n+4}{2}$ is wrong. Hence, we obtain $\an(K_m\Box P_n)>\frac{m(n-1)^2-2n+4}{2}$. That is,  $\an(K_m \Box P_n) \geq \frac{m(n-1)^2-2n+6}{2}$. Thus, the proof is complete by combining the result with Theorem \ref{thm:upper}.    

\textsf{Case-4:} $m \geq 3$ and $n \geq 4$ is even.

The proof follows from Theorem \ref{thm:upper} and \ref{thm:lower}. \qed
\end{proof}

\begin{example}\label{exa1} In Table \ref{tab:K10P7}, a vertex ordering $z_1,z_2,\ldots,z_{70}$ of $V(K_{10} \Box P_7)$, and an optimal radio antipodal labeling of $K_{10} \Box P_7$ are shown.
\begin{table}[hbt!]
\centering
\caption{$\an(K_{10} \Box P_7) = $176.}\label{tab:K10P7}
\begin{tabular}{|c|l|l|l|l|l|l|l|l|l|l|l|l|}
\hline
$(x_i,y'_j)$ & \multicolumn{1} {|c|}{$y'_1$} & \multicolumn{1} {|c|}{$y'_2$} & \multicolumn{1} {|c|}{$y'_3$} & \multicolumn{1} {|c|}{$y'_4$} & \multicolumn{1} {|c|}{$y'_5$} & \multicolumn{1} {|c|}{$y'_6$} & \multicolumn{1} {|c|}{$y'_7$} \\ \hline
    $x_1$ & $z_{47} (118)$ & $z_{65} (165)$ & $z_{35} (86)$ & $z_{53} (133)$ & $z_{1} (0)$ & $z_{41} (101)$ & $z_{59}(148)$ \\ 
    $x_2$ & $z_{54}(136)$ & $z_{2}(3)$ & $z_{42}(104)$ & $z_{60}(151)$ & $z_{8}(18)$ & $z_{48}(119)$ & $z_{66}(166)$ \\
    $x_3$ & $z_{61}(154)$ & $z_{9} (21)$ & $z_{49}(122)$ & $z_{67}(169)$ & $z_{15}(36)$ & $z_{55}(137)$ & $z_{3}(4)$ \\
    $x_4$ & $z_{68}(172)$ & $z_{16}(39)$ & $z_{56}(140)$ & $z_{4}(7)$ & $z_{22}(54)$ & $z_{62}(155)$ & $z_{10}(22)$ \\
    $x_5$ & $z_{5}(10)$ & $z_{23}(57)$ & $z_{63}(158)$ & $z_{11}(25)$ & $z_{29}(72)$ & $z_{69}(173)$ & $z_{17}(40)$ \\
    $x_6$ & $z_{12}(28)$ & $z_{30}(75)$ & $z_{70}(\textbf{176})$ & $z_{18}(43)$ & $z_{36}(90)$ & $z_{6}(11)$ & $z_{24}(58)$ \\ 
    $x_7$ & $z_{19}(46)$ & $z_{37}(93)$ & $z_{7}(14)$ & $z_{25}(61)$ & $z_{43}(108)$ & $z_{13}(29)$ & $z_{31}(76)$ \\
     $x_8$ & $z_{26}(64)$ & $z_{44}(111)$ & $z_{14} (32)$ & $z_{32} (79)$ & $z_{50}(126)$ & $z_{20} (47)$ & $z_{38}(94)$ \\
     $x_9$ & $z_{33}(82)$ & $z_{51}(129)$ & $z_{21} (50)$ & $z_{39} (97)$ & $z_{57} (144)$ & $z_{27}  (65)$ & $z_{45}  (112)$ \\
      $x_{10}$ & $z_{40}   (100)$ & $z_{58}  (147)$ & $z_{28}   (68)$ & $z_{46}  (115)$ & $z_{64}   (162)$ & $z_{34} (83)$ & $z_{52}   (130)$ \\
    \hline
    \end{tabular}
\end{table}
\end{example}
\begin{example}\label{exa2} In Table \ref{tab:K6P8}, a vertex ordering $z_1,z_2,\ldots,z_{48}$ of $V(K_6 \Box P_8)$, and an optimal radio antipodal labeling of $K_6 \Box P_8$ are shown. 
\begin{table}[hbt!]
\small
\centering
\caption{$\an(K_6 \Box P_8) = $150.}\label{tab:K6P8}
\begin{tabular}{|c|l|l|l|l|l|l|l|l|l|l|l|l|}
\hline
$(x_i,y'_j)$ & \multicolumn{1} {|c|}{$y'_1$} & \multicolumn{1} {|c|}{$y'_2$} & \multicolumn{1} {|c|}{$y'_3$} & \multicolumn{1} {|c|}{$y'_4$} & \multicolumn{1} {|c|}{$y'_5$} & \multicolumn{1} {|c|}{$y'_6$} & \multicolumn{1} {|c|}{$y'_7$} & \multicolumn{1} {|c|}{$y'_8$} 
\\ \hline
    $x_1$ & $z_{42}(133)$ & $z_{30}(93)$ & $z_{36}(112)$ & $z_{48}(\textbf{150})$ & $z_{1}(0)$ & $z_{13}(38)$ & $z_{19}(57)$ & $z_{7}(17)$\\ 
    $x_2$ & $z_{2}(3)$ & $z_{14}(41)$ & $z_{20}(60)$ & $z_{8}(20)$ & $z_{25}(78)$ & $z_{37}(116)$ & $z_{43}(135)$ & $z_{31}(95)$\\
    $x_3$ & $z_{26}(81)$ & $z_{38}(119)$ & $z_{44}(138)$ & $z_{32}(98)$ & $z_{9}(26)$ & $z_{21}(64)$ & $z_{3}(5)$ & $z_{15}(43)$\\
    $x_4$ & $z_{10}(29)$ & $z_{22}(67)$ & $z_{4} (8)$ & $z_{16}(46)$ & $z_{33}(104)$ & $z_{45}(142)$ & $z_{27} (83)$ & $z_{39} (121)$\\
    $x_5$ & $z_{34}(107)$ & $z_{46}(145)$ & $z_{28}(86)$ & $z_{40}(124)$ & $z_{17}(53)$ & $z_{5}(12)$ & $z_{11}(31)$ & $z_{23}(69)$\\
    $x_6$ & $z_{18}(55)$ & $z_{6}(15)$ & $z_{12}(34)$ & $z_{24}(72)$ & $z_{41}(130)$ & $z_{29}(90)$ & $z_{35}(109)$ & $z_{47}(147)$\\ 
    \hline
    \end{tabular}
\end{table}
\end{example}


\begin{thebibliography}{00}
\bibitem{Chartrand} G. Chartrand, D. Erwin and P. Zhang, A graph labelling suggested by FM channel restrictions, \emph{Bull. Inst. Combin. Appl.}, \textbf{43} (2005), 43--57.

\bibitem{Zhang1} G. Chartrand, D. Erwin and P. Zhang, Radio antipodal colorings of cycles, \emph{Congr. Numer.}, \textbf{144} (2000), 129--141.

\bibitem{Zhang2} G. Chartrand, D. Erwin and P. Zhang, Radio antipodal colorings of graphs, \emph{Mathematica Bohemica}, \textbf{127} (2002), 57--69.

\bibitem{Zhang3} G. Chartrand, L. Nebeský and P. Zhang, Radio $k$-colorings of paths. \emph{Discuss. Math. Graph
Theory}, \textbf{24} (2004), 5--21.

\bibitem{Griggs} J. R. Griggs and R. K. Yeh, Labelling Graphs with a Condition at Distance $2$, \emph{Discrete Mathematics}, SIAM J. Discrete Math., \textbf{5} (1992), 586--595.

\bibitem{Hale} W. K. Hale, Frequency assignment: Theory and applications, \emph{Proc. IEEE}, \textbf{68} (1980), 1497--1514.

\bibitem{Juan} J. Juan and D. Liu, Antipodal labelings for cycles, \emph{Ars Combinatoria}, \textbf{103} (2012), 81--96.

\bibitem{Togni} R. Khennoufa and O. Togni, A note on radio antipodal colourings of paths, \emph{Mathematica Bohemica}, \textbf{130} (2005), 277--282.

\bibitem{Khennoufa} R. Khennoufa and O. Togni, The radio antipodal and radio numbers of the hypercube, \emph{Ars Combinatoria}, \textbf{102} (2011), 447--461.

\bibitem{Kim} B. M. Kim, W. Hwang, and B. C. Song, Radio number for the product of a path and a complete graph, \emph{J. Comb. Optim.}, \textbf{30} (2015), 139--149.
    
\bibitem{Kumar} K. Kumar and P. Panigrahi, Antipodal number of Cartesian products of complete graphs with cycles. \emph{Communications in Combinatorics and Optimization}, \textbf{10(1)} (2025), 219--232. 

\bibitem{West} D. B. West, Introduction to graph theory, Prentice Hall, 2001.

\bibitem{William} A. William, C. R. Kenneth, Radio antipodal number of certain graphs, \emph{J. Comp. Math. Sci.}, \textbf{2(6)} (2011), 868--872.

\end{thebibliography}
\end{document}